\documentclass[francais]{cedram-ambp}
\usepackage[autolanguage]{numprint}
\usepackage{amsmath,amsfonts,amssymb}
\usepackage{ragged2e}
\usepackage{natbib}
\usepackage[english,frenchb]{babel}
\usepackage[latin1]{inputenc}
\usepackage[T1]{fontenc}
 \usepackage{lmodern}
\usepackage{amsfonts}
\usepackage{amsmath}

\title
{La Z\^eta de Riemann est irrationnelle aux impairs positifs}

\alttitle{Riemann Zeta Function is irrational on Positive odd}

\author
[m. \lastname{kabongo}]
{\firstname{mundankulu} \lastname{kabongo}}

\address{Laboratoire de Recherche Lamunda\\
LouiseVille / Qu\'ebec / Canada}
\email{mundankulu.muvu.kabongo@uqtr.ca}
\altkeywords{Zeta, Riemann, Irrationality, Independence, Linear, Integers, Positives, Values,  Annals, \LaTeX\ class}

\begin{document}
%  Le r{\'e}sum{\'e}
\begin{abstract}

D'abord, nous avons trouv\'{e} une nouvelle \'{e}quation fonctionnelle de la fonction Z\^{e}ta de Riemann. En consid\'erant cette \'equation pour $s=2$ et $s=1$,  nous avons obtenu une relation matricielle qui nous a permis d'obtenir l'irrationalit\'e de $\zeta(j)$ avec $j=3, 4, 5, 6, 7, ...$
\end{abstract}

% Le r�sum� en anglais
\begin{altabstract}
We found a new fonctional equation for Riemann Zeta Function. Considering this equation for $s=2$ and $s=1$, we determinated that for $j=3, 4, 5, 6, 7, ...$; $\zeta(j)$ is irrational.
\end{altabstract}

%% Attention, les r�sum�s sont plac�s *avant* \maketitle
\maketitle

%  D{\'e}but du texte de l'article
\section{Introduction}

Selon la litt\'{e}rature sur la fonction Z\^{e}ta de Riemann, il existe tr\`{e}s peu de r\'{e}sultats sur la nature de cette fonction aux entiers impairs positifs. Alors, nous avons men\'{e} des recherches pour tenter de d\'{e}couvrir quelque chose sur ce sujet. Comme m\'{e}thode, nous avons cherch\'{e} \`{a} trouver des \'{e}quations fonctionnelles en partant d'une g\'{e}n\'{e}ralisation qu'est la fonction Z\^{e}ta d'Hurwitz. Les propri\'{e}t\'{e}s utiles de Z\^{e}ta de Riemann s'y sont d\'{e}coul\'{e}es automatiquement.\\
%  Voici un exemple d'utilisation de l'environnement notation
\begin{notation} 
  Dans cet article, nous d{\'e}signons par  $\zeta(s)$ la fonction Z\^eta de Riemann; $\zeta(s,q)$,  la fonction Z\^eta d'Hurwitz, $\Gamma(s)$ la fonction Gamma d'Euler et $\Longrightarrow$  implique.
\end{notation}
\par
	
	 En effet, de la repr\'{e}sentation en s\'erie de la fonction Z\^{e}ta d'Hurwitz, nous avons trouv\'{e} une autre repr\'esentation, qui nous a permis de construire deux relations contenant les Z\^eta aux entiers positifs. Voici la nouvelle repr\'esentation de la fonction Z\^eta d'Hurwitz.
	
\begin{align}
\zeta(s,q)=\frac{1}{q^{s}}+\sum_{n=0}^{\infty} (-1)^n\frac{q^n}{\Gamma(n+1)}\frac{\Gamma(s+n)\zeta(s+n)}{\Gamma(s)},q\ne 0,-1,-2,...;
\end{align}
$s\ne 1$\\
Dans la partie num\'ero 2 de ce travail nous trouvons la nouvelle \'equation fonctionnelle. La partie num\'ero 3 nous \'ecrivons deux nouvelles combinaisons lin\'eaires ind\'ependantes. Dans la quatri\`eme partie, nous montrons l'irrationalit\'e de z\^eta aux entiers $j=3, 4, 5, 6, 7,...$ Le point 5 est la conclusion.

%Equation fonctionnelle
\section{Nouvelle \'{e}quation fonctionnelle de la fonction Z\^{e}ta de Riemann}
\par
\justify
La fonction Z\^{e}ta d'Hurwitz est donn\'{e}e par la s\'{e}rie:

\begin{equation}
\zeta(s,q)=\sum_{n=0}^\infty \frac{1}{(n+q)^{s}},q\ne 0,-1,-2,-3,.., 
\end{equation}  avec $Re(s)>1$.\\
\\ Trouvons son prolongement analytique.\\
\\Soit la fonction:

\begin{equation}
\Gamma_n(s,q)=\int_0^\infty x^{s-1}e^{-(n+q)x}dx
\end{equation}
Posons                                                                    
$$
x=\frac{u}{n+q}, dx=\frac{du}{n+q}
$$

Cela implique 
                                                                                                                  
$$
\Gamma_n(s,q)=\frac{1}{(n+q)^{s}}\int_0^\infty u^{s-1} e^{-u}du
$$\\
 
  En sommant \`{a} gauche et \`{a} droite on a:                                        
$$
\sum_{n=0}^\infty \Gamma_n(s,q)=\sum_{n=0}^\infty\frac{1}{(n+q)^{s}}\int_0^\infty u^{s-1} e^{-u}du
$$

\begin{equation}
\sum_{n=0}^\infty \Gamma_n(s,q)=\Gamma(s)\zeta(s,q)
\end{equation}\\

Comme ci-dessus, la relation (2.2) est donn\'{e}e par:

$$
\Gamma_n(s,q)=\int_0^\infty x^{s-1} e^{-(n+q)x}dx.
$$
Celle-ci dans (2.3) alors:
                                                                                          
\begin{equation}
\sum_{n=0}^\infty \int_0^\infty x^{s-1} e^{-(n+q)x}dx=\Gamma(s)\zeta(s,q)
\end{equation}

$$
\int_0^\infty e^{-qx} x^{s-1}\sum_{n=0}^\infty e^{-nx}dx=\Gamma(s)\zeta(s,q)
$$                                                                        
$$
\int_0^\infty e^{-qx} x^{s-1} [\frac{1}{(e^x-1)}+1]dx=\Gamma(s)\zeta(s,q)
$$
 Cela implique:                                
\begin{equation}
\int_0^\infty e^{-qx} x^{s-1} \frac{1}{e^x-1}dx+\int_0^\infty e^{-qx} x^{s-1}dx=\Gamma(s)\zeta(s,q).
\end{equation}\\
\\Consid\'erons les $q$ positifs et montrons que la fonction Z\^eta d'Hurwitz converge dans le domaine de la fonction Z\^eta de Riemann.\\
\\En effet, comme:

\begin{equation}
\Bigg|\frac{x^{s-1}}{e^x-1}e^{-qx}\Bigg|\leq\Bigg|\frac{x^{s-1}}{e^x-1}\Bigg|
\end{equation}\\ quelque soit $x\in [0,\infty];$ on a alors Z\^eta d'Hurwitz qui converge dans le domaine de Z\^eta de Riemann.\\
\\Posons et \'{e}crivons
\begin{equation}                                            
x=\frac{u}{q}, dx=\frac{du}{q}, e^{-qx}=\sum_{n=0}^\infty (-1)^n\frac{(q)^n (x)^n}{\Gamma(n+1)}.
\end{equation}

(2.7) dans (2.5) alors:

$$
\sum_{n=0}^\infty (-1)^n\frac{q^n}{\Gamma(n+1)}\int_0^\infty \frac{x^{s-1+n}}{(e^x-1)}dx+\frac{1}{q^{s}}\int_0^\infty u^{s-1} e^{-u}du=\Gamma(s)\zeta(s,q).
$$ 
L'\'equation ci-dessus est \'egale \`a l'\'equation (2.5), par le fait que l'int\'egrale d'une somme est \'egale \` a la somme des integrales.\\
\\Apr\`{e}s int\'{e}gration on trouve:

\begin{equation}
%\fbox{$
\zeta(s,q)=\frac{1}{q^{s}}+\sum_{n=0}^{\infty} (-1)^n\frac{q^n}{\Gamma(n+1)}\frac{\Gamma(s+n)\zeta(s+n)}{\Gamma(s)},s\ne 1;
%$}
\end{equation}

O\`{u} $\Gamma$ est la fonction Gamma d'Euler; $\zeta(s+n)$, la fonction z\^eta de Riemann et $\zeta(s,q)$ la fonction Z\^{e}ta d'Hurwitz.\\
\\(2.8) est une autre repr\'esentation de la fonction Z\^eta d'Hurwitz sous forme d'une s\'erie.\\
\\Selon la relation (2.1), si $q=1$, la fonction Z\^{e}ta d'Hurwitz $\zeta(s,q)$ devient la fonction Z\^{e}ta de Riemann $\zeta(s)$. La relation (2.8) devient donc:

\begin{equation}
%\fbox{$
\zeta(s)=1+\sum_{n=0}^{\infty} (-1)^n\frac{1}{\Gamma(n+1)}\frac{\Gamma(s+n)\zeta(s+n)}{\Gamma(s)}, s\ne 1.
%$}
\end{equation}\\

Dans la relation ci-dessus, le premier terme dans la somme; c'est-\`{a}-dire le terme correspondant \`{a} $n=0$, est \'{e}gale \`{a}:
  
\begin{equation}
\zeta(s)
\end{equation}

(2.10) dans (2.9), on a:

$$
\zeta(s)=1+\zeta(s)+\sum_{n=1}^\infty\frac{(-1)^n}{\Gamma(n+1)}\frac{\Gamma(s+n)\zeta(s+n)}{\Gamma(s)}
$$

Pour $s\ne 1$, cette derni\`{e}re \'{e}quation peut \^{e}tre \'{e}crite comme suite:

\begin{equation}
%\fbox{$
1=\sum_{n=1}^\infty \frac{(-1)^{n+1}}{\Gamma(n+1)}\frac{\Gamma(s+n)\zeta(s+n)}{\Gamma(s)},s\ne 1.
%$}
\end{equation}\\
La relation ci-dessus est la nouvelle \'equation fonctionnelle de la fonction z\^eta de Riemann.

%Combinaisons lineaires

\section{Les relations entre les valeurs de Z\^eta aux entiers positifs; convergence de la s\'erie $K(s)$ pour $s=1$ }

Soient les expressions suivantes:
\begin{equation}
\sum_{n=1}^\infty \frac{(-1)^{n+1}}{\Gamma(n+1)}\frac{\Gamma(s+n)\zeta(s+n)}{\Gamma(s)}=1, s\in C, s \ne 1
\end{equation}

Et 

\begin{equation}
K(s)=\sum_{n=1}^\infty \frac{(-1)^{n+1}}{\Gamma(n+1)}\frac{\Gamma(s+n)\zeta(s+n)}{\Gamma(s)}, s\in C
\end{equation}\\
\\Si on remplace dans la premi\`ere, $s=2$; et dans la deuxi\`eme, $s=1$; on trouve les expressions respectives

\begin{equation}
\sum_{m=3}^\infty (-1)^{m-1}(m-1)\zeta(m)=1
\end{equation}

Et 
 
\begin{equation}
K(1)=\sum_{m=2}^\infty (-1)^{m}\zeta(m)
\end{equation}\\
On peut montrer facilement que la quantit\'e K(1) est \'egale \`a  $$\sum_{n=1}^{\infty}\Bigg(\zeta(2n)-\zeta(2n+1)\Bigg)$$\\
\\
Montrons maintenant qu'elle est finie.\\
\\
\\
\\Soit alors  la fonction $F(s)$ donn\'ee par
\begin{equation}
F(s)=\int_{0}^{\infty}\frac{d}{dx}\Big(\frac{x^{s-1}}{e^x-1}\Big)dx
\end{equation}\\
\\
\\
\\On peut montrer que $F(s)$ est nulle pour les $s$ dont les r\'eels sont sup\'erieurs \`a deux; cela en \'ecrivant:
\\
\begin{equation}
F(s)=\lim\limits_{x\rightarrow \infty} \frac{x^{s-1}}{e^x-1}-\lim\limits_{x\rightarrow 0} \frac{x^{s-1}}{e^x-1}
\end{equation}

Nous disons donc

\begin{equation}
0=\int_{0}^{\infty}\frac{d}{dx}\Big(\frac{x^{s-1}}{e^x-1}\Big)dx
\end{equation}
\\
\\Apr\`es d\'erivation on trouve:

\begin{equation}
0=(s-1)\int_{0}^{\infty}\frac{x^{s-2}}{e^x-1}dx-\int_{0}^{\infty}\frac{x^{s-1}e^x}{(e^x-1)^2}dx
\end{equation}
\\
\\Cela implique que:
\begin{equation}
(s-1)\int_{0}^{\infty}\frac{x^{s-2}}{e^x-1}dx=\int_{0}^{\infty}\frac{x^{s-1}e^x}{(e^x-1)^2}dx
\end{equation}
\\
\begin{equation}
(s-1)\int_{0}^{\infty}\frac{x^{s-2}}{e^x-1}dx=\sum_{m=1}^{\infty}\int_{0}^{\infty}\frac{x^{s-1}e^x}{(e^x-1)}e^{-mx}dx
\end{equation}
\\
\begin{equation}
(s-1)\int_{0}^{\infty}\frac{x^{s-2}}{e^x-1}dx=\sum_{m=0}^{\infty}\int_{0}^{\infty}\frac{x^{s-1}}{(e^x-1)}e^{-mx}dx
\end{equation}
\\
\begin{equation}
(s-1)\int_{0}^{\infty}\frac{x^{s-2}}{e^x-1}dx=\int_{0}^{\infty}\frac{x^{s-1}}{(e^x-1)}dx+\sum_{m=1}^{\infty}\int_{0}^{\infty}\frac{x^{s-1}}{(e^x-1)}e^{-mx}dx
\end{equation}
\\
\begin{equation}
\Gamma(s)\zeta(s-1)=\Gamma(s)\zeta(s)+\sum_{m=1}^{\infty}\int_{0}^{\infty}\frac{x^{s-1}}{(e^x-1)}e^{-mx}dx
\end{equation}
\\
\begin{equation}
\zeta(s-1)=\zeta(s)+\sum_{m=1}^{\infty}\frac{1}{\Gamma(s)}\int_{0}^{\infty}\frac{x^{s-1}}{(e^x-1)}e^{-mx}dx
\end{equation}
\\
\begin{equation}
\zeta(s-1)=\zeta(s)+\sum_{m=1}^{\infty}\Bigg(\zeta(s,m)-\frac{1}{m^{s}}\Bigg)
\end{equation}
\\
\begin{equation}
\zeta(s-1)=\zeta(s)+\sum_{m=2}^{\infty}\zeta(s,m)
\end{equation}
\\
L'\'equation fonctionnelle ci-dessus peut se g\'en\'eraliser. On a donc:
\\
\begin{equation}
\zeta(s)=\zeta(s+1)+\sum_{m=2}^{\infty}\zeta(s+1,m)
\end{equation}
\\
\begin{equation}
\zeta(s+1)=\zeta(s+2)+\sum_{m=2}^{\infty}\zeta(s+2,m)
\end{equation}
\\
\begin{equation}
................................................................
\end{equation}
\\
\begin{equation}
\zeta(s-1+p)=\zeta(s+p)+\sum_{m=2}^{\infty}\zeta(s+p,m)
\end{equation}
\\
L'addition membre \`a  membre de (3.17), (3.18), (3.19) et (3.20) puis \'elimination, on trouve:
\\
\begin{equation}
\zeta(s)=\zeta(s+p)+\sum_{q=1}^{p}\sum_{m=2}^{\infty}\zeta(s+q,m)
\end{equation}
\\
Si $p=\infty$ alors

\begin{equation}
\zeta(s)=1+\sum_{q=1}^{\infty}\sum_{m=2}^{\infty}\zeta(s+q,m)
\end{equation}
R\'eecrivons (3.17) et (3.22):

\begin{equation}
\zeta(s)=\zeta(s+1)+\sum_{m=2}^{\infty}\zeta(s+1,m)
\end{equation}

\begin{equation}
\zeta(s)=1+\sum_{q=1}^{\infty}\sum_{m=2}^{\infty}\zeta(s+q,m)
\end{equation}\\
\\
\\
Si $s=2n$ dans l'\'equation (3.23) et $s=2t$ dans (3.24), alors:

\begin{equation}
\zeta(2n)=\zeta(2n+1)+\sum_{m=2}^{\infty}\zeta(2n+1,m)
\end{equation}
\\
\begin{equation}
\zeta(2t)=1+\sum_{q=1}^{\infty}\sum_{m=2}^{\infty}\zeta(2t+q,m)
\end{equation}\\
\\
\begin{equation}
\zeta(2n)-\zeta(2n+1)=\sum_{m=2}^{\infty}\zeta(2n+1,m)
\end{equation}
\\
\begin{equation}
\zeta(2t)=1+\sum_{q=1}^{\infty}\sum_{m=2}^{\infty}\zeta(2t+q,m)
\end{equation}\\
En sommant sur (3.27), les \'equations deviennent:
\begin{equation}
K(1)=\sum_{n=1}^{\infty}\Bigg(\zeta(2n)-\zeta(2n+1)\Bigg)=\sum_{n=1}^{\infty}\Bigg(\sum_{m=2}^{\infty}\zeta(2n+1,m)\Bigg)
\end{equation}
\\
\begin{equation}
\zeta(2t)=1+\sum_{q=1}^{\infty}\sum_{m=2}^{\infty}\zeta(2t+q,m)
\end{equation}\\
\\
\\
\\
\\On peut encore \'ecrire:

\begin{equation}
K(1)=\sum_{n=1}^{\infty}\Bigg(\zeta(2n)-\zeta(2n+1)\Bigg)=\sum_{n=1}^{t-1}\Bigg(\sum_{m=2}^{\infty}\zeta(2n+1,m)\Bigg)+\sum_{n=t}^{\infty}\Bigg(\sum_{m=2}^{\infty}\zeta(2n+1,m)\Bigg)
\end{equation}

\begin{equation}
\zeta(2t)=1+\sum_{n=t}^{\infty}\sum_{m=2}^{\infty}\zeta(2n+1,m)+\sum_{n=t}^{\infty}\sum_{m=2}^{\infty}\zeta(2n+2,m)
\end{equation}\\
\\
Cela implique que:
\begin{equation}
K(1)=\sum_{n=1}^{\infty}\Bigg(\zeta(2n)-\zeta(2n+1)\Bigg)=\sum_{n=1}^{t-1}\Bigg(\sum_{m=2}^{\infty}\zeta(2n+1,m)\Bigg)+A
\end{equation}

\begin{equation}
\zeta(2t)=1+A+B
\end{equation}\\
\\
Avec 

\begin{equation}
A=\sum_{n=t}^{\infty}\sum_{m=2}^{\infty}\zeta(2n+1,m)
\end{equation}\\
Et
\begin{equation}
B=\sum_{n=t}^{\infty}\sum_{m=2}^{\infty}\zeta(2n+2,m)
\end{equation}\\
(3.34) montre que la quantit\'e $A$ est finie.\\
\\
\\
\\Comme $A$, $t-1$ et $\zeta(2n)-\zeta(2n+1)=\sum_{m=2}^{\infty}\zeta(2n+1,m)$ sont finies, alors (3.33) montre que $K(1)$ est finie. Ceci compl\`ete la d\'emonstration.

%Etude ++++++++++++++++++++++++++++++++++++++++++++++
\section{\'{E}tude d'irrationalit\'{e} de Z\^eta aux entiers $ 3, 4, 5, 6, 7, ..$}

(3.3) et (3.4) peuvent aussi s'\'ecrire par\\
\\
\begin{equation}
a_j\zeta(j)+P_j=1
\end{equation}

\begin{equation}
b_j\zeta(j)+Q_j=-\zeta(2)
\end{equation}

Avec 

\begin{equation}
a_j= (-1)^{j-1}(j-1)
\end{equation}

\begin{equation}
b_j= (-1)^{j}
\end{equation}

\begin{equation}
P_j=\sum_{m=3, m\ne j}^\infty (-1)^{m-1}(m-1)\zeta(m)
\end{equation}

\begin{equation}
Q_j=\sum_{m=3, m\ne j}^\infty (-1)^{m}\zeta(m)-K(1)
\end{equation}
\\
Sous forme matricielle on aura:

\[ \Bigg( \begin{array}{cc}
a_j & P_j \\
\\b_j & Q_j \end{array} \Bigg)
 \Bigg( \begin{array}{c}
\zeta(j) \\
\\1 \end{array}\Bigg) =
\Bigg( \begin{array}{c}
1\\-\zeta(2) \end{array} \Bigg)\]\\
 \\On sait que la famille $\Bigg(1,-\zeta(2)\Bigg)$ est lin\'eairement ind\'ependante sur les rationnels.\\
\\ 
Donc $a-b\zeta(2)=0\Longrightarrow a=b=0$\\
Donc $a\Bigg(a_j\zeta(j)+P_j\Bigg)+b\Bigg(b_j\zeta(j)+Q_j\Bigg)=0\Longrightarrow a=b=0$\\
Donc $\Bigg(a.a_j+b.b_j\Bigg)\zeta(j)+a.P_j+b.Q_j=0\Longrightarrow a=b=0$\\
Donc  la famille $\zeta(j), P_j, Q_j$ est lin\'eairement ind\'ependante sur les rationnels.\\
Donc $P_j, Q_j$ est lin\'eairement ind\'ependante sur les rationnels.\\
\\A partir de l'\'equation matricielle ci-dessus, on peut \'ecrire $P_j$ et $Q_j$\\

\[ \Bigg( \begin{array}{cc}
-a_j &1 \\
\\-b_j &-\zeta(2) \end{array} \Bigg)
 \Bigg( \begin{array}{c}
\zeta(j) \\
\\1 \end{array}\Bigg) =
\Bigg( \begin{array}{c}
P_j\\Q_j \end{array} \Bigg)\]\\
Comme $P_j$, $Q_j$ est lin\'eairement ind\'ependante, alors $p.P_j+q.Q_j=0\Longrightarrow p=q=0$\\
Donc  $p\Bigg(1-a_j\zeta(j)\Bigg)+q\Bigg(-\zeta(2)-b_j\zeta(j)\Bigg)=0\Longrightarrow p=q=0$\\
Donc $-\Bigg(p.a_j+q.b_j\Bigg)\zeta(j)+p-q\zeta(2)=0\Longrightarrow p=q=0$\\
On conclut que $\zeta(j),1,-\zeta(2)$ est une famille lin\'eairement ind\'ependante sur les rationnels.\\
\\
Finalement pour $j=3,4,5,6,7,...$, les $\zeta(j)$ sont tous irrationnels. Par cons\'equent la fonction Z\^eta de Riemann prend des valeurs irrationnelles sur les entiers positifs impairs.

\section{Conclusion}

\justify

La nature de la fonction Z\^{e}ta de Riemann a \'{e}t\'{e} \'{e}tudi\'{e}e de fa\c{c}on s\'{e}par\'{e}e, pour les entiers pairs et impairs. En ce qui nous concerne, nous avons condens\'{e} notre \'{e}tude au moyen d'une seule s\'erie $K(s)$, des fonctions $K_{n}(s)$, qui nous a permis de connaitre l'irrationalit\'{e}, de toutes les valeurs de la fonction Z\^eta aux entiers positifs $j=3, 4, 5, 6, 7, ...$.Ceci en supposant $\zeta(2)$ irrationnel. Ces r\'{e}sultats confirment ceux qui ont \'{e}t\'{e} obtenus avant nous pour les entiers positifs pairs et montrent  l'Irrationalit\'e de Z\^eta aux entiers positifs Impairs.\\ 
\section{R\'ef\'erences}
Nous nous sommes bas\'es sur les d\'efinitions fondamentales. Nous avons pr\'ef\'er\'e l'approche des d\'emonstrations et non l'approche de citations des r\'ef\'erences.

\end{document}